\documentclass{amsart}
\usepackage{amssymb}
\begin{document}
\newtheorem{theorem}{Theorem}[section]
\newtheorem{lemma}[theorem]{Lemma}
\newtheorem{definition}[theorem]{Definition}
\newtheorem{corollary}[theorem]{Corollary}
\def\operatorname#1{{\rm#1\,}}
\def\cupd{\buildrel\ldotp\over\sqcup}
\def\qedbox{\hbox{$\rlap{$\sqcap$}\sqcup$}}
\def\Pspec{\operatorname{spec}}
\def\range{\operatorname{sange}}
\def\Pspan{\operatorname{span}}
\def\RR{{\mathcal{R}}}
\def\RRR{{\text{\pbglie R}}}
\def\rank{\operatorname{rank}}
\def\id{\operatorname{Id}}
\def\trace{\operatorname{trace}}
\def\imag{\operatorname{Im}}
\font\pbglie=eufm10
\def\SS{{\text{\pbglie S}}}
\def\spec{\operatorname{Spec}}
 \makeatletter
  \renewcommand{\theequation}{%
   \thesection.\arabic{equation}}
  \@addtoreset{equation}{section}
 \makeatother
\def\pix{\rho_{\mathcal{X}}}
\title[Spectral geometry of the Curvature Tensor]{The spectral geometry of the\\ Riemann curvature tensor}
\author[P. Gilkey, R. Ivanova, and T. Zhang]{Peter B. Gilkey${}^1$, Raina
Ivanova${}^2$, and Tan Zhang${}^2$}
\thanks{${}^1$Research partially supported by the NSF 
(USA) and MPI (Leipzig)}
\thanks{${}^2$Research supported by the  NSF (USA)}
\begin{address}{PG: Mathematics Department, University of Oregon, Eugene OR 97403 USA.\newline\phantom{...........} email:
gilkey@darkwing.uoregon.edu}\end{address}
\begin{address}{RI: Mathmatics Department,
University of Hawaii - Hilo,
200 W. Kawili St.,\newline\phantom{..........}
Hilo  HI 96720 USA\newline\phantom{..........} email: rivanova@hawaii.edu}\end{address}
\begin{address}{TZ: Department of Mathematics and Statistics,
Murray State University, Murray\newline\phantom{...........}KY 42071 USA\newline\phantom{..........}
email: tan.zhang@murraystate.edu}\end{address}
\date{12 June 2002 Version 1-w SGRCT}
\begin{abstract} Let $\mathcal{E}$ be a natural operator associated to the curvature tensor of a pseudo-Riemannian manifold. We study when the
spectrum, or more generally the real Jordan normal form, of $\mathcal{E}$ is constant on the natural domain of definition. In particular, we examine
the Jacobi operator, the higher order Jacobi operator, the Szabo operator, and the skew-symmetric curvature operator.
\end{abstract}
\maketitle
\font\pbglie=eufm10
\def\Gr{\text{Gr}}

\section{Introduction}\label{Sect1} Let $(M,g)$ be a pseudo-Riemannian manifold of signature $(p,q)$ and dimension $m=p+q$. The associated Riemann
curvature tensor $R$ has the symmetries:
\begin{equation}\begin{array}{l}
R(x,y,z,w)=R(z,w,x,y)=-R(y,x,z,w),\\
R(x,y,z,w)+R(y,z,x,w)+R(z,x,y,w)=0.\end{array}\label{Eqn1.1}\end{equation}
The covariant derivative $\nabla R$ has the symmetries:
\begin{equation}\begin{array}{l}
\nabla R(x,y,z,w;v)=\nabla R(z,w,x,y;v)=-\nabla R(y,x,z,w;v),\\
\nabla R(x,y,z,w;v)+\nabla R(y,z,x,w;v)+\nabla R(z,x,y,w;v)=0,\\
\nabla R(x,y,z,w;v)+\nabla R(x,y,w,v;z)+\nabla R(x,y,v,z;w)=0.\end{array}\label{Eqn1.2}\end{equation}
Instead of working in the geometric category, one can also work in a purely algebraic category. If $V$ is a vector space of signature $(p,q)$, then
we say that $R\in\otimes^4V^*$ is an {\it algebraic curvature tensor} on $V$ if $R$ satisfies the
symmetries of equation (\ref{Eqn1.1}) and that $\RRR\in\otimes^5V^*$ is an {\it algebraic covariant derivative curvature tensor} on $V$ if $\RRR$
satisfies the symmetries of equation (\ref{Eqn1.2}). 

Such tensors are, in general, very complicated objects to study. There are several endomorphisms
$\mathcal{E}$ one can associate to $R$ or to $\nabla R$. One studies $\Pspec(\mathcal{E})$ (i.e. the set of complex eigenvalues of $\mathcal{E}$)
or, more generally, the Jordan normal form of $\mathcal{E}$. In this brief note, we survey some recent
results for the Jacobi operator, the higher order
Jacobi operator, the skew-symmetric curvature operator (in both the real and complex settings), and the Szab\'o operator.

We work in the algebraic setting for the moment. Let
$$S^\pm(V):=\{v\in V:(v,v)=\pm1\}$$
be the {\it pseudo-spheres} of unit  spacelike ($+$) and timelike ($-$) vectors in a vector space $V$ of signature $(p,q)$. Let $\Gr_k(V)$ (resp.
$\Gr_k^+(V)$) be the {\it Grassmannian} of all unoriented (resp. oriented) non-degenerate
$k$ dimensional subspaces of $V$;
$$\Gr_k(V)=\cupd_{r+s=k}\Gr_{r,s}(V)\quad\text{ and }\quad\Gr_k^+(V)=\cupd_{r+s=k}\Gr_{r,s}^+(V)$$
may be decomposed as a disjoint union, where $\Gr_{r,s}(V)$ (resp. $\Gr_{r,s}^+(V)$) is the set of unoriented (resp. oriented)
non-degenerate subspaces of signature
$(r,s)$. We say that a pair $(r,s)$ is {\it admissible} if $\Gr_{r,s}(V)$ consists of more than a single point, or, equivalently, if $0\le r\le p$,
$0\le s\le q$, and $1\le r+s\le m-1$. We say $k$ is admissible if $1\le k\le m-1$.

\begin{definition}\label{Defn1.1}\rm Let $R$ be an algebraic curvature tensor and let $\RRR$ be an algebraic covariant derivative curvature
tensor on a vector space $V$ of signature $(p,q)$. Let $R(\cdot,\cdot)$ and $\RRR(\cdot,\cdot,\cdot)$ be the associated endomorphisms of $V$ defined
by the identities:
$$\begin{array}{l}
(R(x_1,x_2)x_3,x_4)=R(x_1,x_2,x_3,x_4),\quad\text{ and }\quad\\
(\RRR_{x_1}(x_2,x_3)x_4,x_5)=\RRR(x_2,x_3,x_4,x_5;x_1).\end{array}$$
\begin{enumerate}
\smallbreak\item The {\it Jacobi operator} $J(x):y\rightarrow R(y,x)x$ is a self-adjoint map of $V$.  One says that $R$ is
{\it spacelike Osserman} or {\it timelike Osserman} if $\Pspec(J)$ is constant on $S^+(V)$ or on $S^-(V)$, respectively. Similarly, one
says that $R$ is {\it spacelike Jordan Osserman} or {\it timelike Jordan Osserman} if the Jordan normal form of $J$ is constant on $S^+(V)$ or on
$S^-(V)$, respectively.
\smallbreak\item Let $\{e_1,...,e_k\}$ be an orthonormal basis for $\pi\in\Gr_k(V)$. The {\it higher order Jacobi operator}
$J(\pi):=\textstyle\sum_{1\le i\le k}(e_i,e_i)J(e_i)$ is a self-adjoint map of $V$ which is independent of the particular orthonormal basis
chosen for $\pi$. One says that $R$ is {\it Osserman of type
$(r,s)$} if $\Pspec(J)$ is constant on
$\Gr_{r,s}(V)$. The notion {\it Jordan Osserman of type $(r,s)$} is defined similarly.
\smallbreak\item Let $\{e_1,e_2\}$ be an oriented orthonormal basis for $\pi\in\Gr_2^+(V)$. The {\it skew-symmetric curvature operator}
$\RR(\pi):=R(e_1,e_2)$ is independent of the particular orthonormal basis chosen for $\pi$. One says that $R$ is {\it spacelike IP}, {\it mixed
IP}, or {\it timelike IP} if $\Pspec(\RR)$ is constant on $\Gr_{0,2}^+(V)$, $\Gr_{1,1}^+(V)$, or
$\Gr_{2,0}^+(V)$, respectively. The notions {\it spacelike Jordan IP}, {\it timelike Jordan IP}, and {\it mixed Jordan IP} are defined similarly.
\smallbreak\item The {\it Szab\'o operator} $\SS(x):y\rightarrow\RRR_x(y,x)x$ is self-adjoint.
 One says that $\RRR$ is {\it spacelike Szab\'o} or {\it timelike Szab\'o} if $\Pspec(\SS)$ is constant on $S^+(V)$ or on $S^-(V)$,
respectively; the notions {\it spacelike Jordan Szab\'o} and {\it timelike Jordan Szab\'o} are defined similarly.
\smallbreak\item The eigenvalue $0$ plays a distinguished role in this subject. We say that $R$ is {\it 2-nilpotent} if
$R(x_1,x_2)R(x_3,x_4)=0$ for any $x_i\in V$; this implies that $R$ is Osserman, Osserman of type $(r,s)$, and IP.  Similarly, $\RR$ is said to be
{\it 2-nilpotent} if $\RR(x_1,x_2,x_3)\RR(x_3,x_4,x_5)=0$ for any $x_i\in V$; this implies that $\RR$ is Szab\'o.
\end{enumerate}
\end{definition}

The names Osserman, Szab\'o, and IP are used because the germinal work for
this subject in the Riemannian ($p=0$) setting was done by Osserman \cite{refOss}, by Szab\'o
\cite{refSzabo}, and by Ivanov and Petrova \cite{refIP}. 

In the geometric context, we shall say that a pseudo-Riemannian manifold $(M,g)$ has a given property
{\it pointwise} if the Riemann curvature tensor $R$ or the covariant derivative $\nabla R$ has this property on the tangent space $T_PM$ for every
point
$P\in M$. Thus, for example, to say that $(M,g)$ is pointwise spacelike Jordan IP is to mean that $R$ is spacelike Jordan IP on
$T_PM$ for every point
$P\in M$. In particular, the eigenvalues and the Jordan normal form are allowed to vary from point to point. We omit the qualifier `pointwise' if
the structures are independent of the point in question - i.e. if the spectrum or the Jordan normal form does not vary from point to point of the
manifold. This is a subtle, but crucial, distinction.

The following result is proved using analytic continuation - see \cite{refGilkey-01-b} for details:

\begin{theorem}\label{Thm1.2} Let $R$ be an algebraic curvature tensor and let $\RRR$ be an algebraic covariant derivative curvature
tensor on a vector space $V$ of signature $(p,q)$.
\begin{enumerate}
\item The following conditions are equivalent and if any is satisfied, then $R$ is said to be an Osserman tensor:
\begin{enumerate}
\item If $p\ge1$, then $R$ is timelike Osserman;
\item If $q\ge1$, then $R$ is spacelike Osserman.\end{enumerate}
\item The following conditions are equivalent and if any is satisfied, then $R$ is said to be a $k$-Osserman tensor:
\begin{enumerate}
\item There exists $(r,s)$ with $r+s=k$ so $R$ is Osserman of type $(r,s)$;
\item $R$ is Osserman of type $(r,s)$ for any admissible $(r,s)$ with $r+s=k$.
\end{enumerate}
\item The following conditions are equivalent and if any is satisfied, then $R$ is said to be an IP tensor:
\begin{enumerate}
\item If $p\ge2$, then $R$ is timelike IP;
\item If $p\ge1$ and if $q\ge1$, then $R$ is mixed IP;
\item If $q\ge2$, then $R$ is spacelike IP.
\end{enumerate}
\item The following conditions are equivalent and if any is satisfied, then $\RRR$ is said to be a Szab\'o tensor:
\begin{enumerate}
\item If $p\ge1$, then $R$ is timelike Szab\'o;
\item If $q\ge1$, then $R$ is spacelike Szab\'o.
\end{enumerate}
\end{enumerate}
\end{theorem}

By Theorem \ref{Thm1.2}, timelike Osserman and spacelike Osserman are equivalent conditions. Since the Jordan normal form determines the eigenvalue
structure, spacelike Jordan Osserman or timelike Jordan Osserman implies Osserman. The implications are not reversible. As we
shall see in Section
\ref{sect2}, there exist tensors which are Osserman but which are neither spacelike Jordan Osserman nor timelike Jordan Osserman. There also exist
tensors which are spacelike Jordan Osserman but not timelike Jordan Osserman and tensors which are timelike Jordan Osserman but not spacelike
Jordan Osserman. Similar observations hold for the other operators; in the higher signature setting, the eigenvalue structure does not determine
the Jordan normal form of a self-adjoint or of a skew-adjoint linear map.

The following family of examples will play a central role in our discussion. We refer to \cite{refGilkeyIvanovaZhang-02-a} for
details. 
\begin{definition}\label{Defn1.3} Let $u\ge2$ and let $(x,y):=(x_1,...,x_u,y_1,...,y_u)$ be the usual coordinates on $\mathbb{R}^{2u}$. Let
$\psi=\psi_{ij}(x)$ be a symmetric $2$ tensor on
$\mathbb{R}^u$. We define a metric $g_\psi$ of neutral signature $(u,u)$ on $\mathbb{R}^{2u}$ by setting:
$$g_\psi:=\textstyle\sum_{1\le i\le u}dx^i\circ dy^i+\sum_{1\le i,j\le u}\psi_{ij}(x)dx^i\circ dx^j.$$
Let $\mathbb{R}^{(a,b)}$ be Euclidean space with the canonical flat metric $g_{a,b}$ of signature $(a,b)$.
Give $\mathbb{R}^{2u}\times\mathbb{R}^{(a,b)}$ the product metric of signature $(p,q)=(u+a,u+b)$:
$$g_{\psi,a,b}:=g_\psi+g_{a,b}\quad\text{ on }\quad\mathbb{R}^{2u}\times\mathbb{R}^{(a,b)}.$$
\end{definition}

\begin{definition}\label{Defn1.4} If $f$ is a smooth function on $\mathbb{R}^u$, let $\psi_f:=df\circ df$ and let 
$g_{f,a,b}:=g_{\psi_f,a,b}$.
The metric $g_{f,a,b}$ is geometrically realizable as a hypersurface in $\mathbb{R}^{(u+a,u+1+b)}$ with second fundamental form given by the
Hessian $H_{ij}(f):=\partial_i^x\partial_j^xf$.
\end{definition}

\begin{theorem}\label{Thm1.4} Let $g_{\psi,a,b}$ be the metric of Definition \ref{Defn1.3}. Then:
\begin{enumerate}
\item $g_{\psi,a,b}$ is Ricci flat, 2-nilpotent, and Einstein; 
\item $g_{\psi,a,b}$ is neither locally homogeneous nor locally symmetric for generic $\psi$;
\item $g_{\psi,a,b}$ is Osserman, $k$ Osserman for any admissible $k$, IP, and Szab\'o.
\end{enumerate}
\end{theorem}

We will examine these manifolds in greater detail in subsequent sections. We shall almost always assume $p\le q$ so the spacelike directions dominate; the case $q\le p$ can be studied by reversing the sign of the
quadratic form in question. 

Here is a brief guide to this paper. In Section
\ref{sect2}, we study the Jacobi operator and in Section
\ref{sect3}, we study the higher order Jacobi operator. We study the skew-symmetric curvature operator in the real setting in Section
\ref{sect4} and in the complex setting in Section
\ref{sect5}. We conclude in Section
\ref{sect6} by studying the Szab\'o operator.

\section{The Jacobi Operator}\label{sect2}

The classification of Osserman manifolds is essentially complete in the Riemannian ($p=0$) and the Lorentzian ($p=1$) settings. We refer to
\cite{refChia,refNik} for the proof of assertion (1) and to \cite{refBokanBlazicGilkey,refGra} for the proof of assertion (2) in the following
Theorem:

\begin{theorem}\label{Thm2.1}\
\begin{enumerate}
\item Let $(M,g)$ be a Riemannian Osserman manifold of dimension $m\ne16$. Then $(M,g)$ is either flat or locally isometric to a rank $1$
symmetric space.
\item Let $(M,g)$ be a Lorentzian Osserman manifold. Then $(M,g)$ has constant sectional curvature.
\end{enumerate}
\end{theorem}

If $(M,g)$ is either flat or locally isometric to a rank $1$ symmetric space, then $(M,g)$ is both spacelike
Jordan Osserman and timelike Jordan Osserman. Thus one is interested in finding  examples of pseudo-Riemannian manifolds
which are either spacelike Jordan Osserman or timelike Jordan Osserman, but which are neither flat nor local rank $1$ symmetric
spaces.

\medbreak If $W$ is an auxiliary vector space and if $A$
is a linear map of $V$, then we define the {\it stabilization}
$$A\oplus 0:=\left(\begin{array}{ll}A&0\\0&0\end{array}\right)\quad\text{ on }\quad V\oplus W.$$
The Jordan normal form of a spacelike Jordan
Osserman algebraic curvature can be arbitrarily complicated in the balanced setting \cite{refGilkey-01-b}:

\begin{theorem}\label{Thm2.2}
Let $\mathcal{J}$ be a $r\times r$ real matrix and let $p\equiv0$ mod $2^r$. If $V$ is a vector space of neutral
signature
$(p,p)$, then there exists an algebraic curvature tensor $R$ on $V$ so that $J(x)$ is conjugate to
$\pm\mathcal{J}\oplus0$ for every $x\in S^\pm(V)$.
\end{theorem}

The situation is very different if
$p\ne q$. If $p<q$, then the spacelike directions dominate in a certain sense and spacelike Jordan Osserman is a very strong condition. We refer to
\cite{refGilkeyIvanova-02-c} for the proof of the following result:

\begin{theorem}\label{Thm2.3} Let $R$ be an algebraic curvature tensor on a vector space $V$ of signature $(p,q)$. If $p<q$ and if $R$ is spacelike
Jordan Osserman, then $J(x)$ is diagonalizable for every $x\in S^+(V)$.
\end{theorem}

The following result \cite{refGilkeyIvanovaZhang-02-a} deals with the opposite setting and provides examples of nilpotent timelike Jordan Osserman
manifolds if $2\le p\le q$:

\begin{theorem}\label{Thm2.4}  Let $f$ be a smooth function on $\mathbb{R}^u$ with positive definite Hessian. Let $g_{f,a,b}$ be the metric of
signature $(u+a,u+b)$ described in Definition
\ref{Defn1.4}. Then
$g_{f,a,b}$ is timelike Jordan Osserman if and only if $a=0$.
\end{theorem}

Following Garci\'a-Ri\'o, Kupeli, and V\'azquez-Lorenzo
\cite{refGRKVL} (see page 147), we can describe a
family which arises from affine geometry. Let
$\Gamma_{ij}{}^k(x)$ be the Christoffel symbols of an arbitrary torsion free connection
$\nabla$ on $\mathbb{R}^u$. This connection is said to be {\it affine Osserman} if and only if the associated  Jacobi operator is nilpotent.
Define a metric $g_\nabla$ on $\mathbb{R}^{2u}$ by setting:
$$ds^2_\nabla=\textstyle\sum_idx^i\circ dy^i-2\sum_{ijk}y_k\Gamma_{ij}{}^k(x)dx^i\circ dx^j.
$$
 Then
$(M,g_\nabla)$ is Osserman if and only if
$\nabla$ is affine Osserman.
This metric is quite different as the coefficients depend on the $y$ variables as well as on the $x$ variables; there does not
seem to be any direct connection between this metric and the metrics described in Definition \ref{Defn1.3}.
We refer to \cite{refBBGZ,refBCG,refGVV} for other examples of Osserman manifolds in the higher signature setting.

\section{The higher order Jacobi operator}\label{sect3}

There is a basic duality result for the higher order Jacobi operator \cite{refGilkey-01-b,refGSV,refSV}:
\begin{theorem}\label{Thm3.1} Let $R$ be an algebraic curvature tensor on a vector space $V$ of signature $(p,q)$.
\begin{enumerate}\item If $R$ is $k$ Osserman, then $R$ is $m-k$ Osserman.
\item If $R$ is Jordan Osserman of type $(r,s)$, then $R$ is Jordan Osserman of type $(p-r,q-s)$.
\end{enumerate}
\end{theorem}

As $m-1$ Osserman is equivalent to $1$ Osserman, we suppose $2\le k\le m-2$. 
The classification of $k$  Osserman pseudo-Riemannian manifolds is complete in the Riemannian and Lorentzian settings \cite{refGilkey-01-a,refGiSt}:

\begin{theorem} Let $(M,g)$ be a $k$ Osserman pseudo-Riemannian manifold of signature $(p,q)$, where $2\le k\le m-2$. If $p=0$ or if $p=1$, then
$(M,g)$ has constant sectional curvature.
\end{theorem}

Let $\Psi$ be the set of all symmetric $2$ tensors $\psi$ on $\mathbb{R}^u$ so that the Jacobi operator $J_\psi(v)$ defined by the metric $g_\psi$
is positive semi-definite of rank
$p-1$ for every
non zero tangent vector $v$ in $\Pspan\{\partial_i^x\}$. For example, if $f$ is a smooth function on $\mathbb{R}^p$ and if $H(f)$ is positive
definite, then
$df\circ df\in\Psi$. The set $\Psi$ is non-empty, convex, and conelike; it is open in the measurant
topology \cite{refGilkeyIvanovaZhang-02-b}. We have:

\begin{theorem}\label{Thm3.4} Let $\psi\in\Psi$ and let $g_{\psi,a,b}$ be the metric of Definition \ref{Defn1.3}. Then:
\begin{enumerate}
\item  $g_{\psi,a,b}$ is nilpotent Jordan Osserman\begin{enumerate}
\item  of types $(r,0)$ and $( p-r, q)$ if $a=0$ and if $0<r\le p$;
\item of types $(0,s)$ and $( p, q-s)$ if $b=0$ and if $0<s\le p$;
\item of types $(r,0)$ and $( p-r, q)$ if $a>0$ and if  $a+2\le r\le p$;
\item of types $(0,s)$ and $( p, q-s)$ if $b>0$  and if $b+2\le s\le q$;\end{enumerate}
\item $g_{\psi,a,b}$ is not Jordan Osserman of type $(r,s)$ otherwise.
\end{enumerate}
\end{theorem}

See also Bonome et. al. \cite{refBCG} for another family of higher order Osserman manifolds.

\medbreak There is a final family \cite{refGilkeyIvanova-02-a} of algebraic curvature tensors which are higher order nilpotent Jordan Osserman which
do not seem to be realizable geometrically. Let $2\le p\le q$. Let
$\{e_1^-,...,e_p^-,e_1^+,...,e_q^+\}$ be an orthonormal basis for $V$, where the vectors
$\{e_1^-,...,e_p^-\}$ are timelike and the vectors $\{e_1^+,...,e_q^+\}$ are spacelike. 
Let $a$ be a positive integer with $2a\le p$. We define a skew-adjoint linear map
$\Phi_a$ of
$V$ and associated curvature operator $R_a$ by setting:
\begin{equation}\begin{array}{l}\phi_ae_k^\pm=\left\{\begin{array}{lll}
\pm(e_{2i}^-+e_{2i}^+)&\text{if}&k=2i-1\le2a,\\
\mp(e_{2i-1}^-+e_{2i-1}^+)&\text{if}&k=2i\le2a,\\
0&\text{if}&k>2a.\end{array}\right.\\ \\
R_a(x,y)z=(y,\Phi_az)\Phi_ax-(x,\Phi_az)\Phi_ay-2(x,\Phi_ay)\Phi_az.
\end{array}\label{Eqn3.1}\end{equation} 
\begin{theorem}\label{Thm3.5} Let $2a\le p\le q$. Let $R_a$ be the $4$ tensor defined in (\ref{Eqn3.1}). Then:
\begin{enumerate}
\item $R_a$ is an algebraic curvature tensor which is nilpotent $k$ Osserman for any admissible $k$.
\item If $2a<p$, then $R_a$ is Jordan Osserman of type $(p,0)$ and $(0,q)$;
$R_a$ is not Jordan Osserman of type $(r,s)$ otherwise.
\item If $2a=p<q$, then $R_a$ is Jordan Osserman of type $(r,0)$ and of type
$(r,q)$ for any $1\le r\le p-1$; $R_a$ is not Jordan Osserman otherwise.
\item If $2a=p=q$, then $R_a$ is Jordan Osserman of type $(r,0)$, of type
$(r,q)$, of type $(0,s)$, and of type $(p,s)$ for $1\le r\le p-1$ and $1\le s\le q-1$; $R_a$ is not Jordan
Osserman otherwise.
\end{enumerate}
\end{theorem}

\section{The skew-symmetric curvature operator}\label{sect4}

Bounding the rank of a spacelike Jordan IP algebraic curvature tensor is crucial. We refer to
\cite{refGLS} for the proof of assertion (1a), to \cite{refZhang2002} for the proof of assertion (1b), to  \cite{refSt} for the proof of
assertion (1c), and to \cite{refGilkeyZhang2002-a} for the proof of assertions (2) and (3) in the following result.

\begin{theorem} Let $R\ne0$ be an algebraic curvature tensor on a vector space of signature $(p,q)$, where $q\ge5$. 
\begin{enumerate}
\item If
$\RR$ has constant rank $r$ on $\Gr_{0,2}^+(V)$, then
\begin{enumerate}
\item If $p=0$ and if $q=5,6$ or if $q\ge9$, then $r=2$;
\item If $p=1$, and if $q\ge9$, then $r=2$;
\item If $2\le p\le\frac{q-6}4$ and if $\{q,q+1,...,q+p\}$ does not contain a power of $2$, then $r=2$.
\end{enumerate}
\item If $\RR$ has constant rank $2$ on $\Gr_{0,2}^+(V)$, then there exists a self-adjoint map $\phi$ of $V$ whose kernel contains no
spacelike vectors so that $R=\pm R_\phi$, where
$$R_\phi(x,y,z,w):=(\phi x,w)(\phi y,z)-(\phi x,z)(\phi y,w).$$
\item If $\RR$ is spacelike Jordan IP of rank $2$, then $R=\pm R_\phi$, where $\phi$ is self-adjoint and where
$\phi^2=C\cdot\operatorname{Id}$.
\end{enumerate}
\end{theorem}

We refer to
\cite{refGilkey-99} for results concerning the exceptional cases $(p,q)=(0,7)$ and $(p,q)=(0,8)$ in the Riemannian setting.

\medbreak Ivanov and Petrova \cite{refIP} classified the IP Riemannian manifolds of dimension $m=4$.
This classification was subsequently generalized by a number of authors. The following result follows in the Riemannian setting from
\cite{refGLS} and in the pseudo-Riemannian setting from \cite{refGilkeyZhang2002-b}. 

\begin{theorem}\label{Thm4.2} Let $(M,g)$ be a spacelike Jordan IP pseudo-Riemannian manifold of signature $(p,q)$, where $q\ge5$. Assume that
$R_g$ is not nilpotent for at least one point of $M$ and that $\RR$ has constant rank $2$. Then either
$(M,g)$ has constant sectional curvature or locally we can express $ds^2=\varepsilon dt^2+f(t)ds_\kappa ^2$, where $ds_\kappa^2$ has a metric of
constant sectional curvature $\kappa$, where
$f(t)=\varepsilon
\kappa t^2+At+B$, where $\varepsilon=\pm1$, and where
$A^2-4\varepsilon\kappa B\ne0$.\end{theorem}

Theorem \ref{Thm4.2} focuses attention on the nilpotent spacelike Jordan Osserman manifolds since the classification is incomplete in this
setting. We refer to
\cite{refGilkeyIvanovaZhang-02-a} for the proof of the following result:

\begin{theorem}\label{Thm4.3}  Let $f$ be a smooth function with non-degenerate Hessian and let $g_{f,a,b}$ be the metric described in Definition
\ref{Defn1.4}. Then $g_{f,a,b}$ is:
\begin{enumerate}
\item nilpotent IP;
\item never mixed Jordan IP;
\item spacelike Jordan IP if and only if $b=0$;
\item timelike Jordan IP if and only if $a=0$.
\end{enumerate}
\end{theorem}

\section{The skew-symmetric curvature operator in the complex setting}\label{sect5}
 We suppose given a Hermitian almost complex structure $J$ on
$V$, i.e. a skew-adjoint endomorphism $J$ of $V$ which is an isometry such that $J^2=-\operatorname{Id}$. We use
$J$ to give a complex structure to $V$ by defining $\sqrt{-1}x:=Jx$. A
$2$ plane $\pi$ is said to be a {\it complex line} if $J\pi=\pi$; let $\mathbb{CP}(V)$ be the space of all non-degenerate complex lines. A linear
transformation $T$ is said to be {\it complex} if $TJ=JT$.
An algebraic curvature tensor
$R$ is said to be {\it almost complex} if
$J\RR(\pi)=\RR(\pi)J$ for every non-degenerate complex line $\pi$; i.e. if $\RR(\cdot)$ is
complex on $\mathbb{CP}(V)$. Equivalently, \cite{refGilkeyIvanova-01-a} $R$ is almost complex if we have the curvature identity:
$$J^*R=R\quad\text{i.e.}\quad R(x,y,z,w)=R(Jx,Jy,Jz,Jw)\quad\forall(x,y,z,w)\in V.$$
If the Jordan normal form of such a tensor
$\RR(\cdot)$ is constant on $\mathbb{CP}(V)$, then $R$ is said
to be {\it almost complex Jordan IP}. 

Even in the Riemannian setting, the classification is incomplete. We refer to \cite{refGilkey-01-c} for the proof of the following result which
controls the eigenvalue structure:

\begin{theorem}\label{arefm} Let $V$ have signature $(0,q)$. Let $R$ be
an almost complex Jordan IP algebraic curvature tensor on $V$. Let $\{\lambda_s,\mu_s\}$ be
the eigenvalues and multiplicities of $J\RR$, where we order the multiplicities $\mu_s$ so that
$\mu_0\ge...\ge\mu_\ell$. Suppose that
$\ell\ge1$. If $q\equiv2$ mod $4$, then $\ell=1$ and $\mu_1=1$. If $q\equiv0$ mod $4$, then
either $\ell=1$ and $\mu_1\le2$ or $\ell=2$ and $\mu_1=\mu_2=1$.
\end{theorem}

This theorem is sharp, there exist almost complex Jordan IP algebraic curvature tensors with the indicated eigenvalue structures, see
\cite{refGilkeyIvanova-01-a} for details. We remark that Kath \cite{refkath} has obtained some partial
results  concerning almost complex Jordan IP algebraic curvature tensors in signatures
$(p,q)=(0,4)$ or $(p,q)=(2,2)$.

\section{The Szab\'o Operator}\label{sect6}

The classification of Szab\'o tensors is complete in the Riemannian and Lorentzian settings; there are no non-trivial examples. We refer to Szab\'o
\cite{refSzabo} for the case $p=0$ and to Gilkey-Stavrov \cite{refGiSt} for case $p=1$ in the  following result:

\begin{theorem}\label{Thm6.1} Let $\RRR$ be a Szab\'o tensor on a vector space of signature $(p,q)$, where $p=0$ or $p=1$. Then
$\RRR=0$.
\end{theorem}

Szab\'o used this observation to give a proof of the well known fact that any $2$ point homogeneous Riemannian manifold is locally symmetric;
similarly any pointwise totally isotropic Lorentzian manifold is locally symmetric.

Theorem \ref{Thm1.4} provides examples of non-trivial Szab\'o tensors if $2\le p\le q$. However, these examples are neither spacelike Jordan
Szab\'o nor timelike Jordan Szab\'o if
$\nabla R\ne0$. In fact, there are no known examples of spacelike Jordan Szab\'o or timelike
Jordan Szab\'o algebraic curvature tensors such that $\RRR\ne0$ and it is natural to conjecture that none exist. There are some partial results in
this direction. Let
$A$ be self-adjoint and let
$\lambda\in\mathbb{C}$. We define real operators
$A_\lambda$ on $V$ and associated generalized subspaces $E_\lambda$ by setting:
\begin{eqnarray}
  &&A_\lambda:=\left\{\begin{array}{ll}
A-\lambda\cdot\id&\quad\text{if}\quad\lambda\in\mathbb{R},\\
(A-\lambda\cdot\id)(A-\bar\lambda\cdot\id)&\quad\text{if}\quad\lambda\in\mathbb{C}-\mathbb{R},\end{array}\right.\nonumber\\
  &&E_\lambda:=\ker\{A_\lambda^m\}.\label{eqn1.1}\end{eqnarray}
As $A_\lambda=A_{\bar\lambda}$, we see $E_\lambda=E_{\bar\lambda}$. Both $A$ and $A_\lambda$ preserve each
generalized eigenspace $E_\lambda$. The operator $A$
is said to be {\it Jordan simple} if $A_\lambda=0$ on $E_\lambda$ for all $\lambda$. 

Let $\nu(q)$ be the {\it Adams number}; this is the maximal
number of linearly independent vector fields on the sphere $S^{q-1}$ in $\mathbb{R}^q$ \cite{refAdams}. If $\RR$ is a Szab\'o tensor, then let
$\Pspec^\pm(\SS)$ denote the spectrum of the associated 

Szab\'o operator on $S^\pm(V)$. We refer to \cite{refGilkeyIvanova-02-b} for the proof of
the following theorem in the algebraic context:

\begin{theorem}\label{Thm6.2} Let $\RRR$ be an algebraic covariant derivative curvature tensor on a vector space of signature $(p,q)$. 
\begin{enumerate}\item If $\RRR$ is Szab\'o, then:\begin{enumerate}
\item $\Pspec^\pm\{\SS\}=-\Pspec^\pm\{\SS\}=\sqrt{-1}\phantom{.}\Pspec^\mp\{\SS\}$;
\item $\Pspec^\pm\{\SS\}\subset\mathbb{R}\cup\sqrt{-1}\phantom{.}\mathbb{R}$;
\item If $p<q$, then $\Pspec^+\{\SS\}\subset\sqrt{-1}\phantom{.}\mathbb{R}$ and $\Pspec^-\{\SS\}\subset\mathbb{R}$.
\end{enumerate}
\item If $p<q$ and if $\RRR$ is spacelike Jordan Szab\'o, then:\begin{enumerate}
\item $\SS(v)$ is Jordan simple for any $v\in S^+(V)$;
\item If $p<q-\nu(q)$, then $\rank\{\SS(v)\}\le2\nu(q)$ for any $v\in S^+(V)$;
\item If $q$ is odd, then $\RRR=0$.
\end{enumerate}
\end{enumerate}
\end{theorem}

We say that a pseudo-Riemannian manifold $(M,g)$ is {\it pointwise totally
isotropic} if given a point $P\in M$ and nonzero tangent vectors $X$ and $Y$ in $T_PM$ with
$(X,X)=(Y,Y)$, there is a local isometry of
$(M,g)$, fixing $P$, which sends $X$ to $Y$. The Szab\'o
operator of a pointwise totally isotropic pseudo-Riemannian manifold necessarily has constant Jordan normal form on $S^+(T_PM)$
for any $P\in M$. Thus Theorem \ref{Thm6.2} yields the following Corollary:

\begin{corollary}\label{cor1.3} Let $(M,g)$ be a pointwise totally isotropic pseudo-Riemannian manifold of
signature $(p,q)$, where $p<q$ and where $q$ is odd. Then
$(M,g)$ is locally symmetric.
\end{corollary}

More generally, Wolf
showed that any locally totally isotropic manifold is necessarily locally symmetric, see \cite{refWolf} (Theorem 12.3.1). 
Thus Corollary \ref{cor1.3} gives a different proof of Wolf's result if $p<q$ and if $q$ is odd.\vfill\eject

\end{document}